\documentclass[11pt]{article}
\usepackage{amssymb, amsmath,  amsfonts,dsfont, mathrsfs,euscript,eufrak}
\usepackage{textcomp}
\usepackage{latexsym}
\usepackage{indentfirst}
\usepackage[normalem]{ulem}

\usepackage{graphicx}
\usepackage[colorlinks=true, allcolors=blue]{hyperref}

\newtheorem{Theorem}{Theorem}
\newtheorem{Fact}{Fact}
\newtheorem{Lemma}{Lemma}

\newtheorem{Corollary}{Corollary}

\newcommand{\dd}{\,d}
\newcommand{\id}{\mathds{1}}
\newcommand{\iid}{\mathbb{I}}

\newcommand{\Complex}{\mathbb C}
\newcommand{\R}{\mathbb R}
\newcommand{\Z}{\mathbb Z}
\newcommand{\N}{\mathbb N}
\newcommand{\Q}{\mathbb Q}

\newcommand{\RealPart}{\mathrm{Re}}
\newcommand{\ImagPart}{\mathrm{Im}}

\newcommand{\DF}{\boldsymbol{F}}
\newcommand{\DFabs}{\boldsymbol{F}_{\!a}}
\newcommand{\DFcont}{\boldsymbol{F}_{\!c}}
\newcommand{\DFdisc}{\boldsymbol{F}_{\!d}}
\newcommand{\DFdiscL}{\boldsymbol{F}_{\!d}^{\mathrm{lat}}}
\newcommand{\DFdiscNL}{\boldsymbol{F}_{\!d}^{\mathrm{nl}}}

\newcommand{\ID}{\boldsymbol{I}}
\newcommand{\RID}{\boldsymbol{Q}}
\newcommand{\RIDabs}{\boldsymbol{Q}_a}
\newcommand{\RIDsing}{\boldsymbol{Q}_s}
\newcommand{\RIDdisc}{\boldsymbol{Q}_d}
\newcommand{\RIDdiscL}{\boldsymbol{Q}_d^{\mathrm{lat}}}
\newcommand{\RIDdiscNL}{\boldsymbol{Q}_d^{\mathrm{nl}}}
\newcommand{\RIDmixdabs}{\boldsymbol{Q}_{d,a}^{\mathrm{mix}}}
\newcommand{\RIDmixlabs}{\boldsymbol{Q}_{\mathrm{lat},a}^{\mathrm{mix}}}

\newcommand{\lext}{\mathrm{lext}}
\newcommand{\rext}{\mathrm{rext}}
\newcommand{\cext}{\mathrm{cext}}
\newcommand{\Supp}{\mathcal{S}}

\newcommand{\e}{\varepsilon}

\newcommand{\U}{\mathcal{U}}

\begin{document}

\title{Denseness in total variation and the class of rational-infinitely divisible distributions}

\author{A. A. Khartov$^{1,2,}$\footnote{Email addresses: \texttt{khartov.a@iitp.ru}, \texttt{alexeykhartov@gmail.com}} }

\footnotetext[1]{Institute for Information Transmission Problems (Kharkevich Institute) of Russian Academy of Sciences, Bolshoy Karetny per. 19, build.1, 127051 Moscow, Russia.}
\footnotetext[2]{Smolensk State University, 4 Przhevalsky st., 214000 Smolensk, Russia. }

\maketitle
\begin{abstract}
	We study a new class of so-called rational-infinitely (or quasi-infinitely) divisible probability laws on the real line. The characteristic functions of these  distributions are ratios of the characteristic functions of classical infinitely divisible laws and they admit L\'evy--Khinchine type representations with ``signed spectral measures''. This class is rather wide and it has a lot of nice properties. For instance, this class is dense in the family of all (univariate) probability laws with respect to weak convergence. In this paper, we consider the questions concerning a denseness of this class with respect to convergence in total variation. The problem is considered separately for different types of probability laws taking into account the supports of the distributions. A series of ``positive'' and ``negative'' results are obtained.	
\end{abstract}

\textit{Keywords and phrases}:  infinite divisibility, rational-infinite divisibility, the L\'evy--Khinchine type representation,  denseness, total variation.

\section{Introduction}
Let us consider the family of all probability laws on the real line. One of the most important classes in this family is the class of infinitely divisible laws (see \cite{GnedKolm}). Their theory is deep and beautiful (the L\'evy--Khinchine formula, subclasses of self-decomposable and stable laws, decomposition problems). Their significance for probability and stochastic processes is fundamental and great (limit theorems for independent random variables, the L\'evy processes, stochastic calculus).  The applications are rather extensive (physics, teletraffic models, insurance mathematics).

Let us recall the definition and the main properties of infinitely divisible laws. Let $F$ be  a (right-continuous) distribution function on the real line $\R$ with the characteristic function 
\begin{eqnarray*}
	f(t):=\int_{\R} e^{itx} \dd F(x),\quad t\in\R.
\end{eqnarray*}
Recall that $F$ (and the corresponding law) is called \textit{infinitely divisible} if for every  $n\in\N$ (set of positive integers) there exists a  distribution function $F_{1/n}$ such that $F=(F_{1/n})^{*n}$, where ``$*$'' denotes a convolution operation, i.e. $F$ is $n$-fold convolution power of $F_{1/n}$. So $f(t)=(f_{1/n}(t))^n$, $t\in\R$, where $f_{1/n}$ is the characteristic function of $F_{1/n}$, $n\in\N$. It is well known that $f(t)\ne 0$ for any $t\in\R$. The most important fundamental fact is that $F$ is infinitely divisible if and only if $f$ is represented by \textit{the L\'evy--Khinchine formula}:
\begin{eqnarray}\label{repr_f}
	f(t)=\exp\biggl\{it \gamma+\int_{\R} \bigl(e^{itx} -1 -it \sin(x)\bigr)\tfrac{1+x^2 }{x^2} \dd G(x)\biggr\},\quad t\in\R,
\end{eqnarray}
with some \textit{shift parameter} $\gamma\in\R$, and with a  bounded non-decreasing \textit{spectral function} $G: \R \to \R$  that is assumed to be right-continuous at every point of the real line with the condition $G(-\infty)=0$.  \textit{The spectral pair} $(\gamma, G)$ is uniquely determined by $f$ and hence by $F$. Using spectral pairs, it is possible to describe the various properties of $F$ and characterize weak convergence of infinitely divisible laws (see \cite{Petr} and \cite{Steutel}). 

Let $\DF$ be the family of all distribution functions on the real line. We will denote by $\ID$ the class of all infinitely divisible distribution functions. Note that the class $\ID$ is rather wide in $\DF$ and it contains important well known distribution functions. Nevertheless, a lot of probability laws are not infinitely divisible, for instance, non-degenerate distributions with bounded support (see \cite{BaxterShapiro}). In addition, sequences of distribution functions from $\ID$ can not approximate (say converge to)  non-infinitely divisible distribution functions even in weak sense, because $\ID$ is closed under weak convergence (see \cite{Sato1999}, Lemma~7.8).  

There is a class of distributions functions, which is very similar to $\ID$ and  it is additionally dense in $\DF$ with respect to weak convergence. Let us give precise definitions. We call a distribution function $F$ (and the corresponding law) \textit{rational-infinitely divisible} if there are some infinitely divisible distribution functions $F_1$ and $F_2$ such that $F_1=F*F_2$. This equality is equivalent to $f(t)=f_1(t)/f_2(t)$, $t\in\R$, where $f_1$ and $f_2$ are characteristic functions of $F_1$ and $F_2$, respectively. It is clear that here  $f(t)\ne 0$ for any $t\in\R$, because $f_1$ does not vanish and $f_2$ is bounded (both on $\R$). It is not difficult to see that $f$ admits \textit{the L\'evy--Khinchine type representation} \eqref{repr_f} with $\gamma=\gamma_1-\gamma_2$ and $G=G_1-G_2$, where $(\gamma_1, G_1)$ and $(\gamma_2, G_2)$ denote the spectral pairs of $F_1$ and $F_2$, respectively. It is important to note that the function $G$ has a bounded variation on $\R$ and it is non-monototic in general.  In addition, $G$ is  right-continuous at every point and $G(-\infty)=0$.  The  pair $(\gamma, G)$ is uniquely determined by $F$ as for the class $\ID$ (it was proved \cite{GnedKolm}, p. 80). Now we inversely suppose that $f$ is represented by  \eqref{repr_f} with some $\gamma\in\R$ and  $G$ satisfying the above-mentioned conditions. Following the terminology from \cite{LindPanSato}, a distribution function $F$ is called  \textit{quasi-infinitely divisible} in this case. Then $F$ will be rational-infinitely divisible by the Hahn--Jordan decomposition for the function $G$ of bounded variation. It should be noted that, due to the L\'evy--Khinchine type representations, class $\RID$ finds interesting applications in various fields (see reviews and references in \cite{BergKutLind}, \cite{KhartovMixSingPart}, and \cite{KhartovNecess}).

We will denote by $\RID$ the class of all rational-(quasi-)infinitely divisible distribution functions. It is clear that $\RID\ne \DF$, because $\RID$ doesn't contain distribution functions, whose characteristic function vanish at some points. However, according to the above, we have $\ID\subset \RID$. There are particular examples of distribution functions from $\RID\setminus\ID$, which appeared even in \cite{GnedKolm} (see p.~81--83)  and  \cite{LinOstr} (see p.~165).  Moreover, $\RID$ is a very significant extension of $\ID$ that is partially confirmed by the facts below. Let us denote by $\DFdisc$ the set of all discrete distribution functions and let $\DFdiscL$  be  the subclass of its lattice representatives. The following very nice result was proved in \cite{LindPanSato} (Theorem~8.1 and Corollary~8.2).
\begin{Fact}\label{fact_discL}
	Suppose that $F\in\DFdiscL$. Then $F\in\RID$ if and only if $f(t)\ne 0$ for any $t\in\R$.
\end{Fact}
In the paper \cite{Berger}, this assertion was extended  to an important class of mixture distributions. Here and below $\DFabs$ denotes the set of all absolutely continuous distribution functions. 
\begin{Fact}\label{fact_discLabs}
	Suppose that $F(x)=c_dF_d(x)+(1-c_d) F_a(x)$, $x\in\R$, where $F_d\in\DFdiscL$, $F_a\in\DFabs$, and $c_d\in(0,1]$. Let $f(t)=c_d f_d(t)+(1-c_d) f_a(t)$, $t\in\R$, where $f_d$ and $f_a$ are the characteristic functions of $F_d$ and $F_a$, respectively. Then $F\in\RID$ if and only if $f_d(t)\ne 0$ and  $f(t)\ne 0$ for any $t\in\R$.
\end{Fact}
There are several recent results, which complement these propositions. For instance, Fact~\ref{fact_discL} was generalized in \cite{AlexeevKhartov} and \cite{KhartovSPL2022}. Since it will be used below, we recall it here. 
\begin{Fact}\label{fact_disc}
	Suppose that $F\in\DFdisc$. Then $F\in\RID$ if and only if $\inf_{t\in\R}|f(t)|>0$.
\end{Fact}
The analog of Fact~\ref{fact_discLabs} for $F_d\in\DFdisc$ can be find in \cite{BergerKutlu}. Mixtures with continuous singular parts were considered in \cite{KhartovMixSingPart}. There is a special approach in \cite{KhartovGeneral}, which is convenient  for pure absolutely continuous distribution functions in some cases. Simpler conditions of belonging to the class $\RID$ and a lot of other important results and examples concerning this class are given in the above-mentioned paper \cite{LindPanSato}, which remains the most important reference on this topic.

The class $\RID$ is dense in $\DF$ with respect to weak convergence, as was proved in \cite{LindPanSato} (see Theorem~4.1). In fact, the authors constructively showed even more, namely, that the subclass $\RIDdiscL:=\RID\cap \DFdiscL$ is dense in $\DF$ in weak sense, i.e. for any $F\in\DF$  we can construct a sequence $(F_n)_{n\in\N}$ from $\RIDdiscL$ (the authors provided a method for it), which weakly converges to $F$.  We note, however, that this denseness property is lost if we turn to the multivariate analog of $\RID$ (see \cite{BergKutLind}) as it follows from one assertion by Kutlu (see \cite{Kutlu}).

From the above, we get at once the following denseness results with respect to (formally stronger) uniform convergence in the class $\DFcont$ of all continuous distribution functions. Indeed,  if $(F_n)_{n\in\N}$ from $\RIDdiscL$ weakly converges to a function $F\in\DFcont$, then this sequence uniformly converges to $F$ by the  well known fact (see \cite{Petr}, Theorem~1.11). Thus $\RIDdiscL$ and hence $\RIDdisc$ are dense in $\DFcont$ in that sense. Next, let $(\sigma_n)_{n\in\N}$ be a sequence of positive reals such that $\sigma_n\to 0$ as $n\to\infty$ and let $\Phi_{\sigma_n}$, $n\in\N$, denote the Gaussian distribution functions with zero means and variances $\sigma_n$, $n\in\N$, respectively. Then the sequence $(F_n*\Phi_{\sigma_n})_{n\in\N}$ uniformly converges to $F$ too, but here all $F_n*\Phi_{\sigma_n}$ are absolutely continuous functions from $\RID$. It means that $\RIDabs:=\RID\cap\DFabs$ is dense in $\DFcont$ with respect to uniform convergence. In addition, the class $\RIDsing$ of all continuous singular representatives of $\RID$ is dense in $\DFcont$ in that sense. In order to show it, we take a continuous singular infinitely divisible distribution function $H_s$ (it is known that such a function exists, see \cite{Tucker}) and define $H_{s,\sigma_n}(x):=H_s(x/\sigma_n)$, $x\in\R$, for every $n\in\N$. Since $(H_{s,\sigma_n})_{n\in\N}$ weakly converges to the degenerate distribution function (with unit jump at $x=0$), the sequence $(F_n*H_{s,\sigma_n})_{n\in\N}$ weakly and hence uniformly convergence to continuous $F$. Here the functions $F_n*H_{s,\sigma_n}$ belong to $\RID$, and they are continuous singular as  convolutions of discrete $F_n$ and continuous singular $H_{s,\sigma_n}$, $n\in\N$. This yields the required assertion. By the way, it seems that these remarks have not been discussed before in papers. 

It is natural and interesting to consider the questions concerning a denseness of the class $\RID$ with respect to convergence in total variation on $\R$ (for simplicity, we will say \textit{denseness in variation}). However,  we are not aware of any general results here and so we want to devote the article to these problems. Namely, we will invesigate an existence of such denseness property  for the most important subclasses $\RIDdisc$ and $\RIDabs$ in the sets $\DFdisc$ and $\DFabs$, correspondingly. Moreover, we will obtain the answer about the denseness in variation of the whole class $\RID$ in the family $\DF$. The results will be rather non-trivial and even surprising. 

There will be two groups of results. The first group consists of ``positive'' results, where the denseness property holds (see Section~3). The results for which it is not true (say ``negative'' results) are assigned to another group (see Section~4).  All these results are obtained due to the tools from Section~2. Some additional notation will be introduced in these sections before formulating the results.

\section{The main tools}

Let us start with  necessary notation. Following \cite{LinOstr}, for any distribution function $F$ we define the set of all points of increase of $F$ (the support of the corresponding distribution)
\begin{eqnarray*}
	\Supp(F):=\bigl\{x\in\R: F(x+\e)-F(x-\e)>0 \text{ for any }\e>0\bigr\}.
\end{eqnarray*}
and its boundaries
\begin{eqnarray*}
	\lext(F):=\inf\Supp(F)\geqslant -\infty,\qquad \rext(F):=\sup\Supp(F)\leqslant \infty.
\end{eqnarray*}
For bounded $\Supp(F)$ we also define its center
\begin{eqnarray*}
	\cext(F):=\dfrac{\lext(F)+\rext(F)}{2}.
\end{eqnarray*}
As usual (see \cite{Lukacs}, p.~30) we call a distribution $F$ \textit{symmetric} if $F(x)=1-F(-x-0)$ for any $x\in\R$, but $F$ is said to be \textit{shift-symmetric} if the distribution function $x\mapsto F(x+c)$, $x\in\R$, is symmetric for some $c\in\R$. In that case, $c=\cext(F)$ for bounded $\Supp(F)$.

We first formulate a useful lemma, which will be a key tool in the proofs of ``positive'' results. In fact, it is a modification of one assertion by Ibragimov from \cite{LinOstr} (see Theorem~7.1.1.). For the convenience of the reader we fully prove this lemma here (besides, the proof in \cite{LinOstr} contains some inaccuracies).

\begin{Lemma}\label{lm_key}
	Let $F_0$ be a  distribution function on $\R$ with a bounded $\Supp(F_0)$.  Let $f_0$ be the characteristic function of  $F_0$.	Let $\gamma_0$ be an arbitrary real number such that $\gamma_0\ne \cext(F_0)$  if $F_0$ is shift-symmetric.
	Then for all $\delta\in(0,1)$ except a countable set (it may be empty, finite or infinite) we have
	\begin{eqnarray}\label{lm_key_conc}
		\delta e^{it\gamma_0}+(1-\delta) f_0(t)\ne 0\quad\text{for any}\quad t\in\R.  
	\end{eqnarray}
	In particular, for any $\tau\in(0,1]$ there exists $\delta\in(0,\tau)$ such that \eqref{lm_key_conc} holds. 
\end{Lemma}
\noindent\textbf{Proof of Lemma~\ref{lm_key}.}
	Let us fix any $\gamma_0$ as in the formulation and define $F_1(x):=F_0(x+\gamma_0)$, $x\in\R$.  The distribution function $F_1$ has the characteristic function $f_1(t):=f_0(t)e^{-it\gamma_0}$, $t\in\R$. The support $\Supp(F_1)$ is bounded and hence $f_1$ is an entire characteristic function (see \cite{LinOstr} p.~35), i.e. there is a function $\tilde{f}_1$, which is analytic in $\Complex$ (set of complex numbers), such that $\tilde{f}_1(t)=f_1(t)$ for any $t\in\R$. Observe that, due to the assumptions on $\gamma_0$, the distribution function $F_1$ is  non-symmetric for any case.  Therefore the function  $t\mapsto\ImagPart f_1(t)$, $t\in\R$, is not an identical zero. Additionaly,  observe that 
\begin{eqnarray*}
	\ImagPart f_1(t)= \dfrac{f_1(t)-\overline{f_1(t)}}{2i}= \dfrac{f_1(t)-f_1(-t)}{2i}=\dfrac{\tilde{f}_1(t)-\tilde{f}_1(-t)}{2i},\quad t\in\R.
\end{eqnarray*}
It is clear that the function $z\mapsto \tfrac{\tilde{f}_1(z)-\tilde{f}_1(-z)}{2i}$, $z\in\Complex$, is an analytic in $\Complex$ and it is not an identical zero. Therefore it has only a countable set of zeroes in $\Complex$. Consequently, $\ImagPart f_1(t)=0$  only for $t$ from some countable set $T\subset\R$ ($T$ may be empty, finite or infinite).

For every $\delta\in(0,1)$ we define the following function
\begin{eqnarray*}
	f_\delta(t) :=\delta+(1-\delta) f_1(t),\quad t\in\R.
\end{eqnarray*}
Let us fix $\delta\in(0,1)$ and observe that 
\begin{eqnarray*}
	f_\delta(t) =\delta+(1-\delta) \RealPart f_1(t)+i(1-\delta)\ImagPart f_1(t),\quad t\in\R.
\end{eqnarray*}
Therefore $f_\delta(t)=0$ if and only if $\delta+(1-\delta) \RealPart f_1(t)=0$ and $\ImagPart f_1(t)=0$. Hence $f_\delta(t)\ne 0$ if $\ImagPart f_1(t)\ne 0$, i.e. $t\notin T$. We now assume that $t'\in T$. So we have  $f_\delta(t') =\delta+(1-\delta) \RealPart f_1(t')$.  If $\RealPart f_1(t')\geqslant 0$, then $f_\delta(t') \geqslant\delta>0$.  If $\RealPart f_1(t')< 0$, then $f_\delta(t') =\RealPart f_1(t')+\delta(1-\RealPart f_1(t'))=0$ only in the case
\begin{eqnarray*}
	\delta=\delta':= \dfrac{-\RealPart f_1(t')}{1-\RealPart f_1(t')}\in(0,1).
\end{eqnarray*}
Let $D$ denote the set of all values of such $\delta'$ for varying $t'\in T$. Thus if $\delta\notin D$ then $f_\delta(t)\ne 0$ and hence $f_\delta(t)e^{it\gamma_0}=\delta e^{it\gamma_0}+(1-\delta) f_0(t)\ne 0$ for any $t\in\R$, i.e. \eqref{lm_key_conc} is valid. Since $T$ is countable, $D$ is a countable set too (it may be empty, finite or infinite). In particular, for any $\tau\in(0,1]$ the set $(0,\tau)\setminus D$ is always uncountable and, in this set,  there exists a $\delta$ for which \eqref{lm_key_conc} holds  as required.\quad $\Box$\\

The next lemmata are used in the proofs of ``negative'' results.  The first is an interesting assertion proved by Kutlu in the short paper \cite{Kutlu} (see Theorem~1).
\begin{Lemma}\label{lm_Kutlu}
	The function
	\begin{eqnarray*}
		\varphi(t_1,t_2):=\dfrac{1}{3}\,e^{it_1}+\dfrac{1}{3}\,e^{it_2}+\dfrac{1}{3}\,e^{i(t_1+t_2)},\quad t_1,t_2\in\R,
	\end{eqnarray*}
	cannot be approximated arbitrarily well by zero-free continuous
	functions from $\R^2$ to $\Complex$ with respect to uniform convergence on $[-\pi,\pi]^2$.
\end{Lemma}

The next lemma is a recent non-trivial result from \cite{KhartovNecess} (see Theorem~2). It complements Facts~\ref{fact_discL}--\ref{fact_disc} from the introduction.
\begin{Lemma}\label{lm_fdsep0}
	Suppose that $F(x)=c_dF_d(x)+(1-c_d) F_c(x)$, $x\in\R$, where $c_d\in[0,1]$, $F_d$ and $F_c$ are discrete and continuous  distribution functions, respectively. Let $f_d$ be the characteristic function of $F_d$.	If $F\in\RID$ and $c_d>0$, then $\inf_{t\in\R}|f_d(t)|>0$.
\end{Lemma}
The proofs of these lemmata can be found in the given references.

\section{The ``positive'' results}
Recall that $\DFabs$ denotes the family of all absolutely continuous distribution functions and $\RIDabs$ denotes the class of rational-infinitely divisible representatives of $\DFabs$, i.e. $\RIDabs=\RID\cap\DFabs$. We additionally introduce their subclasses, which are corresponded to distributions with supports contained in a particular set. Namely, for any $S\subset \R$ we define
\begin{eqnarray*}
	\DFabs(S):= \bigl\{F\in\DFabs:  \Supp(F)\subset S\bigr\}\quad\text{and} \quad\RIDabs(S):= \bigl\{F\in\RIDabs:  \Supp(F)\subset S\bigr\}.
\end{eqnarray*}
In particular, $\DFabs(\R)=\DFabs$ and $\RIDabs(\R)=\RIDabs$. It is also clear that $\RIDabs(S)=\RID\cap\DFabs(S)$.

For any $S\subset\R$ and $\tau>0$ we need closed $\tau$-neighbourhoods of $S$. More precisely, we define
\begin{eqnarray}\label{def_Staupm}
	[S]_{\tau}^+:=\overline{S+[0,\tau]}\quad\text{and} \quad [S]_{\tau}^-=\overline{S+[-\tau,0]}.
\end{eqnarray}
Here $\overline{A}$ means the closure of a set $A\subset \R$ (in the standard topology of $\R$), and the sums of sets are understood in the arithmetic sense, i.e. $A+B:=\{x+y: x\in A, y\in B\}$ for any $A,B\subset\R$. 

In the proofs of the results, we will use the following indicator functions:
\begin{eqnarray*}
	\iid_A(x):=
	\begin{cases}
		1,& x\in A,\\
		0,& x\notin A,
	\end{cases}
	\quad\text{and}\quad
	\id_a(x):=
	\begin{cases}
		1,& x\leqslant a,\\
		0,& x>a,
	\end{cases}
\end{eqnarray*}
with $A\subset \R$ and $a\in\R$. Next, for any function $H$ defined on $\R$ we denote by $\|H\|$ the total variation of $H$ on $\R$.

Now we are ready to formulate and prove the first ``positive'' result.
\begin{Theorem}\label{th_AbsTV}
	Let $S$ be a  subset of $\R$ such that $\DFabs(S)\ne \varnothing$. 
	Then for any $\tau>0$  the classes $\RIDabs([S]_\tau^+)$ and $\RIDabs([S]_\tau^-)$ are dense in variation in the class $\DFabs(S)$.
\end{Theorem}
This theorem obviously implies the following important assertion.
\begin{Corollary}\label{co_th_AbsTV}
	The class $\RIDabs$ is dense in variation in the class $\DFabs$.
\end{Corollary}
\noindent\textbf{Proof of Theorem~\ref{th_AbsTV}.}
Let $F$ be a distribution function from $\DFabs(S)$ with the density $p$. Let us take the density function $b_q(x):=C_q\cdot q^x(1-q)^{1-x}\iid_{[0,1]}(x)$, $x\in\R$, of the continuous Bernoulli distribution (see \cite{Korkmaz}), where $q\in(0,1)$ is the main parameter and $C_q$ is a constant such that $C_q\cdot\int_0^1 q^x(1-q)^{1-x} \dd x=1$. Let $B_q$ denote the corresponding distribution function. Note that we come to the uniform distribution on $[0,1]$ if $q=\tfrac{1}{2}$. But we fix $q\ne \tfrac{1}{2}$ and for any $\tau>0$ we define transformed versions of $B_q$:
\begin{eqnarray*}
	B_{q,\tau}^+(x):=B_q(x/\tau)\quad\text{and}\quad   B_{q,\tau}^-(x):=B_q((x+\tau)/\tau),\quad x\in\R,
\end{eqnarray*}
which have the densities
\begin{eqnarray*}
	b_{q,\tau}^+(x):=\dfrac{b_q(x/\tau)}{\tau}\quad\text{and}\quad   b_{q,\tau}^-(x):=\dfrac{b_q((x+\tau)/\tau)}{\tau},\quad x\in\R,
\end{eqnarray*}
respectively.  It is easily seen that these distributions weakly converge to the degenerate law as $\tau\to 0+$, i.e.  we have
\begin{eqnarray}\label{conv_Bpmqtau}
	B^{\pm}_{q,\tau}(x)\to\id_0(x)\quad \text{as}\quad\tau\to 0\!+\quad\text{for any}\quad x\ne 0.
\end{eqnarray}

We first show that  
\begin{eqnarray}\label{conv_FPhisigma}
	\lim\limits_{\tau \to 0+} \|F-F*B^{\pm}_{q,\tau}\|=0.
\end{eqnarray}
For every $\tau>0$ the distribution functions $F*B_{q,\tau}^{\pm}$ have the densities
\begin{eqnarray*}
	(p*b^{\pm}_{q,\tau})(x)= \int_{\R} p(x-u)b^{\pm}_{q,\tau}(u)\dd u,\quad x\in\R,
\end{eqnarray*}
respectively, and hence we have the formula
\begin{eqnarray*}
	\|F-F*B^{\pm}_{q,\tau}\|=\int_{\R} \bigl|p(x)-(p*b^{\pm}_{q,\tau})(x)\bigr| \dd x.
\end{eqnarray*}
We observe that
\begin{eqnarray*}
	p(x)-(p*b^{\pm}_{q,\tau})(x)&=&p(x)-\int_{\R} p(x-u)b^{\pm}_{q,\tau}(u)\dd u\\
	&=&\int_{\R} \bigl(p(x)-p(x-u)\bigr)b^{\pm}_{q,\tau}(u)\dd u,\quad x\in\R,
\end{eqnarray*}
and, consequently,
\begin{eqnarray*}
	\bigl|p(x)-(p*b^{\pm}_{q,\tau})(x)\bigr|\leqslant \int_{\R} \bigl|p(x)-p(x-u)\bigr| b^{\pm}_{q,\tau}(u)\dd u,\quad x\in\R.
\end{eqnarray*}
Thus for any $\tau>0$ we get
\begin{eqnarray}\label{ineq_varFFPhisigma}
	\bigl\|F-F*B^{\pm}_{q,\tau}\bigr\|&\leqslant &\int_{\R} \biggl(  \int_{\R} \bigl|p(x)-p(x-u)\bigr|b^{\pm}_{q,\tau}(u)\dd u \biggr) \dd x\nonumber\\
	&=&\int_{\R} \Delta(u) b^{\pm}_{q,\tau}(u)\dd u=\int_{\R} \Delta(u) \dd B^{\pm}_{q,\tau}(u),
\end{eqnarray}
where
\begin{eqnarray*}
	\Delta(u):= \int_{\R} |p(x)-p(x-u)|\dd x,\quad u\in\R.
\end{eqnarray*}
Since $p\in L_1(\R)$, it is known that $\Delta(u)\to 0$ as $u\to 0$  (see \cite{MakarovPodkor}, p.~518, and \cite{Rudin}, p.~182). So  $\Delta$ is continuous on $\R$, because for any real $u$ and $u_0$
\begin{eqnarray*}
	|\Delta(u)-\Delta(u_0)|&\leqslant& \int_{\R} \bigl||p(x)-p(x-u)|- |p(x)-p(x-u_0)| \bigr|\dd x\\
	&\leqslant& \int_{\R} |p(x-u_0)-p(x-u)|\dd x\\
	&=& \int_{\R} |p(y)-p(y-(u-u_0))|\dd y=\Delta(u-u_0).
\end{eqnarray*}
Since the function $p$ is a density, $\Delta$ is bounded on $\R$:
\begin{eqnarray*}
	0\leqslant \Delta(u)\leqslant\int_{\R} (p(x)+p(x-u))\dd x=\int_{\R} p(x)\dd x+\int_{\R} p(x-u)\dd x=2,\quad u\in\R.
\end{eqnarray*}
According to these properties of $\Delta$ and the convergence from  \eqref{conv_Bpmqtau}, we obtain
\begin{eqnarray*}
	\lim\limits_{\tau\to 0+} \int_{\R} \Delta(u) \dd B^{\pm}_{q,\tau}(u)=\int_{\R} \Delta(u) \dd\id_0(u),
\end{eqnarray*}  
by the extended Helly's Second Theorem (see \cite{Lukacs} p. 47). The integral is equal to $\Delta(0)=0$ in the right-hand side. Therefore \eqref{conv_FPhisigma} holds due to \eqref{ineq_varFFPhisigma}. Note that \eqref{conv_FPhisigma} may be proved by direct applying of Theorem 9.3.3 from \cite{MakarovPodkor}, but we preferred a more transparent proof.

We now fix $\tau>0$ and $\e\in(0,1)$. According to \eqref{conv_FPhisigma}, there is  $\tau_\e\in(0,\tau)$ such that $\bigl\|F-F*B^{\pm}_{q,\tau_\e}\bigr\|<\e$. Next, we choose $r_\e>0$ such that
\begin{eqnarray}\label{ineq_int_p_e}
	\int_{\R\setminus[-r_\e,r_\e]} p(x)\dd x<\e.
\end{eqnarray}
Let us define $c_\e:=\int_{[-r_\e,r_\e]} p(x)\dd x> 1-\e>0$, and   $p_\e(x):= \tfrac{1}{c_\e}\,p(x)\iid_{[-r_\e,r_\e]}(x)$, $x\in\R$.  It is clear that $p_\e$ is a density of some  distribution function on $\R$, say $F_\e$. Let $f_\e$ denote its characteristic function. It is seen that
$-r_\e\leqslant\lext(F_\e)<\rext(F_\e)\leqslant r_\e$. So $\cext(F_\e)\in[-r_\e,r_\e]$ and we fix any real $\gamma_\e\in \Supp(F)\cap[-r_\e,r_\e]$, but $\gamma_\e\ne \cext(F_\e)$. By Lemma \ref{lm_key}, there exists $\delta_\e\in (0,\e)$ such that $\delta_\e e^{it\gamma_\e}+(1-\delta_\e)f_\e(t)\ne 0$ for any $t\in\R$. The function $f_\e^\circ(t):=\delta_\e e^{it\gamma_\e}+(1-\delta_\e)f_\e(t)$, $t\in\R$, is characteristic for the distribution function  $F_\e^\circ(x):=\delta_\e\id_{\gamma_\e}(x)+(1-\delta_\e) F_\e(x)$, $x\in\R$, which is a mixture of the absolutely continuous function $F_\e$ and the degenerate component $\id_{\gamma_\e}$. By the way, we note that $\Supp(F_\e^\circ)\subset \Supp(F)\subset S$. According to Fact~\ref{fact_discLabs} (where we set  $F_d:=\id_{\gamma_\e}$, $F_a:=F_\e$, and $c_d:=\delta_\e$), $F_\e^\circ\in\RID$ if and only if $e^{it\gamma_\e}\ne 0$ and $f_\e^\circ(t)\ne 0$ for any $t\in\R$. Since the last inequalities are true, it follows that $F_\e^\circ\in\RID$. Next, in the paper \cite{KhartovGeneral} (see Example~3.3, p.~7), it was shown that $B_q\in\RID$ (we fixed $q\ne \tfrac{1}{2}$). Consequently, $B^{\pm}_{q,\tau_\e}\in\RID$ and hence $F_\e^\circ * B^{\pm}_{q,\tau_\e}\in\RID$ too, because the class $\RID$ is closed under  shifting, scaling and convolution (see introduction). The distribution functions $F_\e^\circ * B^{\pm}_{q,\tau_\e}$ are absolutely continuous as a result of convolution with absolutely continuous $B^{\pm}_{q,\tau_\e}$. Thus we have $F_\e^\circ * B^{\pm}_{q,\tau_\e}\in\RIDabs$. Next, we observe (using Theorem~3.2.1. from \cite{LinOstr}) that
\begin{eqnarray*}
	\Supp(F_\e^\circ * B^{\pm}_{q,\tau_\e})=\overline{\Supp(F_\e^\circ)+\Supp(B^{\pm}_{q,\tau_\e})}\subset \overline{S+\Supp(B^{\pm}_{q,\tau})}=[S]_{\tau}^\pm
\end{eqnarray*}
due to $\tau_\e<\tau$, $\Supp(B^{+}_{q,\tau})=[0,\tau]$, $\Supp(B^{-}_{q,\tau})=[-\tau,0]$, and according to the definitions from \eqref{def_Staupm}.  So we conclude that $F_\e^\circ * B^{\pm}_{q,\tau_\e}\in\RIDabs([S]_{\tau}^\pm)$, respectively.

We now estimate the total variation of the difference $F$ and $F_\e^\circ * B^{\pm}_{q,\tau_\e}$. We first observe that
\begin{eqnarray*}
	\bigl\|F- F_\e^\circ * B^{\pm}_{q,\tau_\e}\bigr\|\leqslant \bigl\|F- F * B^{\pm}_{q,\tau_\e}\bigr\|+\bigl\|F * B^{\pm}_{q,\tau_\e}- F_\e^\circ * B^{\pm}_{q,\tau_\e}\bigr\|.
\end{eqnarray*}
Here the first term of the sum is strictly less than $\e$ (see above), and for the second term we have
\begin{eqnarray*}
	\bigl\|F * B^{\pm}_{q,\tau_\e}- F_\e^\circ * B^{\pm}_{q,\tau_\e}\bigr\|=\bigl\|(F- F_\e^\circ)* B^{\pm}_{q,\tau_\e}\bigr\|\leqslant \|F- F_\e^\circ\|\cdot \bigl\| B^{\pm}_{q,\tau_\e}\bigr\|= \|F- F_\e^\circ\|,
\end{eqnarray*}
where we used well-known properties of the convolution  and the total variation (see \cite{GelRaiShil}, p. 165--166). We now consider $F- F_\e^\circ= F-\delta_\e\id_{\gamma_\e}-(1-\delta_\e)F_\e$. It has the discrete part $\delta_\e\id_{\gamma_\e}$ (with the total variation $\delta_\e$) and the absolutely continuous part $F-(1-\delta_\e)F_\e$. Therefore
\begin{eqnarray*}
	\|F-F_\e^\circ\|=\|\delta_\e\id_{\gamma_\e}\|+\|F-(1-\delta_\e) F_\e\|=\delta_\e+\|F-(1-\delta_\e) F_\e\|.
\end{eqnarray*}
Thus we have
\begin{eqnarray}\label{ineq_diffFFePhi}
	\bigl\|F- F_\e^\circ * B^{\pm}_{q,\tau_\e}\bigr\|< \e+\delta_\e+\|F-(1-\delta_\e) F_\e\|.
\end{eqnarray}
Let us estimate the latter variation:
\begin{eqnarray*}
	\|F-(1-\delta_\e) F_\e\|&=&\|F-c_\e F_\e- (1-c_\e-\delta_\e) F_\e\|\\
	&\leqslant& \|F-c_\e F_\e\|+ |1-c_\e-\delta_\e|\cdot \|F_\e\|,
\end{eqnarray*} 
where $\|F_\e\|=1$ (because $F_\e$ is a distribution function). Let us consider the function $F-c_\e F_\e$. For any $x\in\R$
\begin{eqnarray*}
	F(x)-c_\e F_\e(x)&=& \int_{(-\infty,x]} p(u)\dd u-c_\e \int_{(-\infty,x]} \dfrac{1}{c_\e}\, p(u)\iid_{[-r_\e,r_\e]}(u)\dd u\\
	&=&\int_{(-\infty,x]} \bigl(p(u)- p(u)\iid_{[-r_\e,r_\e]}(u)\bigr)\dd u\\
	&=&\int_{(-\infty,x]}p(u)\iid_{\R\setminus [-r_\e,r_\e]}(u)\dd u.
\end{eqnarray*} 
There is non-negative integrand in the last integral. Then, on account of \eqref{ineq_int_p_e}, we obtain
\begin{eqnarray}\label{ineq_FceFe}
	\|F-c_\e F_\e\|= \int_{\R} p(x)\iid_{\R\setminus [-r_\e,r_\e]}(x) \dd x=\int_{\R\setminus[-r_\e,r_\e]} p(x)\dd x<\e.
\end{eqnarray}
So we have $\|F-(1-\delta_\e) F_\e\|<\e+ |1-c_\e-\delta_\e|$. Using this estimate in \eqref{ineq_diffFFePhi}, we get
\begin{eqnarray*}
	\bigl\|F- F_\e^\circ * B^{\pm}_{q,\tau_\e}\bigr\|< 2\e+\delta_\e+ |1-c_\e-\delta_\e|.
\end{eqnarray*}
Recall that $\delta_\e\in(0,\e)$ and $c_\e\in (1-\e,1]$ (see above). So $1-c_\e\in [0,\e)$ and hence $|1-c_\e-\delta_\e|<\e$. Therefore
\begin{eqnarray*}
	\bigl\|F- F_\e^\circ * B^{\pm}_{q,\tau_\e}\bigr\|< 4\e,
\end{eqnarray*}
where $\e$ can be chosen arbitrarily small and $F_\e^\circ * B^{\pm}_{q,\tau_\e}\in\RIDabs([S]_{\tau}^\pm)$ as we showed above.

Thus for any $F\in\DFabs(S)$ we can find an element from $\RIDabs([S]_{\tau}^+)$ or $\RIDabs([S]_{\tau}^-)$ (with any $\tau>0$), which is arbitrarily close to $F$ in variation. It means that the classes $\RIDabs([S]_{\tau}^\pm)$ are dense in variation in the class $\DFabs(S)$.\quad $\Box$\\

We now turn to the family $\DFdisc$  of all discrete distribution functions and to the class $\RIDdisc=\RID\cap\DFdisc$. Similarly to the case of absolutely continuous distribution functions, for any $S\subset \R$ we  introduce the following classes:
\begin{eqnarray*}
	\DFdisc(S):= \bigl\{F\in\DFdisc:  \Supp(F)\subset S\bigr\}\quad\text{and}\quad \RIDdisc(S):= \bigl\{F\in\RIDdisc:  \Supp(F)\subset S\bigr\}.
\end{eqnarray*}
It is clear that $\DFdisc(\R)=\DFdisc$, $\RIDdisc(\R)=\RIDdisc$, and $\RIDdisc(S)=\RID\cap\DFdisc(S)$. Next, it is important for us to divide $\DFdisc$ into two subclasses $\DFdiscL$ and $\DFdiscNL$ of all discrete lattice and non-lattice distribution functions, respectively (see \cite{Petr}, p.~2). Similarly, we divide the class $\RIDdisc$ into $\RIDdiscL:=\RIDdisc\cap \DFdiscL$ and $\RIDdiscNL:=\RIDdisc\cap \DFdiscNL$. In this notation, if $L$ is a non-empty subset of a lattice $\{a+bk:\,  k\in\Z\}$ with some $a\in\R$ and $b>0$, then $\DFdisc(L)\subset \DFdiscL$ and $\RIDdisc(L)\subset \RIDdiscL$.

The second ``positive'' result is obtained for the discrete lattice case.
\begin{Theorem}\label{th_DiscLTV}
	Suppose that $L\subset\{a+bk:\,  k\in\Z\}$ with some $a\in\R$, $b>0$, and $L\ne \varnothing$. Then $\RIDdisc(L)$ is dense in variation in the class $\DFdisc(L)$.
\end{Theorem}
Hence we have the following obvious conclusion.
\begin{Corollary}\label{co_th_DiscLTV}
	The class $\RIDdiscL$ is dense in variation in the class $\DFdiscL$.
\end{Corollary}
\noindent\textbf{Proof of Theorem~\ref{th_DiscLTV}.}
 Let $F$ be a distribution function from $ \DFdisc(L)$ with characteristic function $f$. Then $f$ is represented by the formula
\begin{eqnarray*}
	f(t)=\sum_{k\in\Z} p_k e^{it(a+bk)},\quad t\in\R,
\end{eqnarray*} 
where each $p_k$ denotes the probability of the point $a+bk$ of the corresponding distribution. So $p_k\geqslant 0$ for any $k\in\Z$, $\sum_{k\in\Z} p_k=1$, and $p_k=0$ if $a+bk\notin L$. Without loss of generality we can assume that $F\in \DFdisc(L)\setminus \RIDdisc(L)$. Then $F$ is non-degenerate, because all degenerate distribution functions belongs to $\RID$ (see introduction). This means that $\Supp(F)$ (and hence $L$) contains at least two distinct points $a+bk_1$ and $a+bk_2$ with integers $k_1$ and $k_2$ (we set $k_1<k_2$).

Let us fix arbitrary $\e\in(0,1)$ and choose integers  $K_{1,\e}$ and $K_{2,\e}$ such that   $a+bK_{1,\e}$ and $a+bK_{2,\e}$ belong to $\Supp(F)$,  $K_{1,\e}\leqslant k_1<k_2\leqslant K_{2,\e}$, 
\begin{eqnarray}\label{ineq_sumpkeK12e}
	q_{1,\e}:=\sum_{k\in\Z:\,k<K_{1,\e}} p_k<\dfrac{\e}{2}\quad\text{and}\quad q_{2,\e}:=\sum_{k\in\Z:\,k>K_{2,\e}} p_k<\dfrac{\e}{2}.
\end{eqnarray}
Next, we define
\begin{eqnarray*}
	\tilde{f}_\e(t):= (p_{K_{1,\e}}+q_{1,\e}) e^{a+b K_{1,\e}}+\!\!\!\! \sum_{\substack{k\in\Z:\\ K_{1,\e}<k<K_{2,\e}}}\!\!\!\! p_k e^{it(a+bk)}+(p_{K_{2,\e}}+q_{2,\e}) e^{a+b K_{2,\e}},\quad t\in\R.
\end{eqnarray*}
If the set $\{k\in\Z: K_{1,\e}<k<K_{2,\e}\}$ is empty, then the corresponding sum is assumed to be zero. It is easily seen that $\tilde{f}_\e$ is the characteristic function of some probability distribution, whose function will be denoted by $\widetilde{F}_\e$. Observe that
\begin{eqnarray*}
	\bigl\|F-\widetilde{F}_\e\bigr\|&=&\sum_{k\in\Z:\,k<K_{1,\e}} |p_k-0| + \sum_{\substack{k\in\Z:\\ K_{1,\e}<k<K_{2,\e}}} |p_k-p_k|+\sum_{k\in\Z:\,k>K_{2,\e}}| p_k-0|\\
	&&{}+ |p_{K_{1,\e}}-(p_{K_{1,\e}}+q_{1,\e})| +|p_{K_{2,\e}}-(p_{K_{2,\e}}+q_{2,\e})|\\
	&=& \sum_{k\in\Z:\,k<K_{1,\e}} p_k + \sum_{k\in\Z:\,k>K_{2,\e}} p_k+q_{1,\e}+q_{2,\e}\\
	&=&2(q_{1,\e}+q_{2,\e}).
\end{eqnarray*}
Therefore $\|F-\widetilde{F}_\e\|<2\e$ due to \eqref{ineq_sumpkeK12e}. 
Next, by definitions of $K_{1,\e}$ and $K_{2,\e}$, we have $\Supp(\widetilde{F}_\e)\subset \Supp(F)\subset L$ and 
\begin{eqnarray*}
	\lext(\widetilde{F}_\e)=a+bK_{1,\e}\leqslant a+bk_1<a+bk_2 \leqslant a+bK_{2,\e}=\rext(\widetilde{F}_\e).
\end{eqnarray*}
From this we conclude that $\cext(\widetilde{F}_\e)$ can not be equal to $\lext(\widetilde{F}_\e)$ or $\rext(\widetilde{F}_\e)$. So we set $\gamma_\e:=\lext(\widetilde{F}_\e)=a+bK_{1,\e}$, and, according to  Lemma~\ref{lm_key}, we choose $\delta_\e\in(0,\e)$  such that
\begin{eqnarray*}
	\tilde{f}_\e^\circ (t):= \delta_\e e^{it\gamma_\e}+(1-\delta_\e) \tilde{f}_\e(t)\ne 0\quad\text{for any}\quad t\in\R.
\end{eqnarray*}
Let $\widetilde{F}_\e^\circ$ denote the distribution function corresponded to $\tilde{f}_\e^\circ$, i.e. $\widetilde{F}_\e^\circ(x)=\delta_\e\id_{\gamma_\e}(x)+(1-\delta_\e) \widetilde{F}_\e(x)$, $x\in\R$. By Fact~\ref{fact_discL}, we know that $\widetilde{F}_\e^\circ\in\RID$.  Moreover, since $\gamma_\e\in \Supp(\widetilde{F}_\e)$ and  $\Supp(\widetilde{F}_\e)\subset L$, it is clear that  $\widetilde{F}_\e^\circ\in\RIDdisc(L)$. Next, 
\begin{eqnarray*}
	\bigl\|\widetilde{F}_\e-\widetilde{F}_\e^\circ\bigr\|=\bigl\|\widetilde{F}_\e-\delta_\e \id_{\gamma_\e}-(1-\delta_\e) \widetilde{F}_\e\bigr\|=\bigl\|\delta_\e(\widetilde{F}_\e-\id_{\gamma_\e})\bigr\|\leqslant2\delta_\e<2\e.
\end{eqnarray*}
Therefore
\begin{eqnarray*}
	\bigl\|F-\widetilde{F}_\e^\circ\bigr\|\leqslant	\bigl\|F-\widetilde{F}_\e\bigr\|+\bigl\|\widetilde{F}_\e-\widetilde{F}_\e^\circ\bigr\|<4\e.
\end{eqnarray*}

Thus for any $F\in \DFdisc(L)$ we found $\widetilde{F}_\e^\circ\in\RIDdisc(L)$, which is arbitrarily close to $F$ in variation.\quad $\Box$\\

We supplement the previous theorems by the following ``positive'' result for the mixtures of discrete lattice and absolutely continuous distribution functions. We previously introduce the necessary notation. For any two classes $\DF_{\!1}$ and $\DF_{\!2}$ of distribution functions we denote by $\DF_{\!1}\oplus \DF_{\!2}$ the family of all distribution functions of the form $F(x)=c F_1(x)+(1-c) F_2(x)$, $x\in\R$, where  $c\in[0,1]$, $F_1\in\DF_{\!1}$ and $F_2\in\DF_{\!2}$.  For instance, $\DFdisc\oplus\DFabs$ denotes the class of all mixtures of discrete and absolutely continuous distribution functions. Next, we introduce the classes
\begin{eqnarray*}
	\RIDmixdabs:=\RID\cap (\DFdisc\oplus\DFabs)\quad\text{and} \quad\RIDmixlabs:=\RID\cap (\DFdiscL\oplus\DFabs).
\end{eqnarray*}
Let $L$ and $S$ be subsets of $\R$ such that $\DFdisc(L)$ and $\DFabs(S)$ are non-empty. We define 
\begin{eqnarray*}
	\RIDmixdabs(L,S):= \RID\cap (\DFdisc(L)\oplus\DFabs(S)).
\end{eqnarray*}
It is important to understand that $\RIDmixdabs(L,S)\ne \RIDdisc(L)\oplus \RIDabs(S)$, however, of course,
\begin{eqnarray}\label{conc_RIDdiscRIDabsRIDmix}
	\RIDdisc(L)\subset\RIDmixdabs(L,S)\quad\text{and}\quad  \RIDabs(S)\subset\RIDmixdabs(L,S).
\end{eqnarray}
If $L$ is a non-empty subset of a lattice, then  $\RIDmixdabs(L,S)\subset\RIDmixlabs$ for any admissible $S$. 

\begin{Theorem}\label{th_Mixdabs}
	Let $S$ be a  subset of $\R$ such that $\DFabs(S)\ne \varnothing$. Let $L\subset\{a+bk:\,  k\in\Z\}$ with some $a\in\R$, $b>0$, and $L\ne \varnothing$. The class $\RIDmixdabs(L,S)$ is dense in variation in the class $\DFdisc(L)\oplus\DFabs(S)$.
\end{Theorem}
From this we immediately get the following assertion.
\begin{Corollary} 
	The class $\RIDmixlabs$ is dense in variation in the class $\DFdiscL\oplus \DFabs$.
\end{Corollary}
\noindent\textbf{Proof of Theorem~\ref{th_Mixdabs}.}
Suppose that $F\in\DFdisc(L)\oplus\DFabs(S)$ and $F$ has the characteristic function $f$. Then
\begin{eqnarray*}
	F(x)=c_d F_d(x)+(1-c_d) F_a(x),\quad x\in\R,
\end{eqnarray*}
where $c_d\in[0,1]$, $F_d\in \DFdisc(L)$ and $F_a\in \DFabs(S)$. Let $f_d$ and $f_a$ be the characteristic functions of $F_d$ and $F_a$, respectively. So we have
\begin{eqnarray*}
	f(t)=c_d f_d(t)+(1-c_d) f_a(t),\quad t\in\R.
\end{eqnarray*}
We note the case $c_d=1$ is covered by Theorem~\ref{th_DiscLTV} according to \eqref{conc_RIDdiscRIDabsRIDmix}. So we next assume that  $c_d\in[0,1)$.

\textbf{1)} Suppose that $L$ contains only one point $\gamma_1:=a+bk_1$ with some $k_1\in\Z$. It means that $\Supp(F_d)=\{\gamma_1\}$, i.e. $F_d=\id_{\gamma_1}$, and  
\begin{eqnarray*}
	F(x)=c_d \id_{\gamma_1}(x)+(1-c_d) F_{a}(x),\quad x\in\R.
\end{eqnarray*}
Let $p$ be the density of $F_a$. For any $\e\in(0,\tfrac{1}{2})$ we define a new density function $p_\e$ exactly as in the proof of Theorem~\ref{th_AbsTV} with the same $r_\e$ and $c_\e$ (see the definitions after \eqref{ineq_int_p_e}). Let $F_{a,\e}$ and $f_{a,\e}$ be the distribution function  and the characteristic function corresponding to   $p_\e$, respectively. Note that $\Supp(F_{a,\e})\subset\Supp(F_a)\subset S$ and  $-r_\e\leqslant\lext(F_{a,\e})<\rext(F_{a,\e})\leqslant r_\e$. Next, we define the distribution function
\begin{eqnarray*}
	F_\e(x):=c_d \id_{\gamma_1}(x)+(1-c_d) F_{a,\e}(x),\quad x\in\R.
\end{eqnarray*}
It is clear that $\Supp(F_\e)$ is bounded. Let $f_\e$ be the characteristic function of $F_\e$, i.e. $f_\e(t):=c_de^{it\gamma_1}+(1-c_d)f_{a,\e}(t)$, $t\in\R$. 

\textbf{1.a)} Assume that $\cext(F_\e)\ne \gamma_1$. Then, according to Lemma~\ref{lm_key}, there exists $\delta_\e\in (0,\e)$ such that $\delta_\e e^{it\gamma_1}+(1-\delta_\e)f_{\e}(t)\ne 0$ for any $t\in\R$. Observe that $f_\e^\circ(t):=\delta_\e e^{it\gamma_1}+(1-\delta_\e)f_\e(t)$, $t\in\R$, is the  characteristic function for the distribution function  $F_\e^\circ(x):=\delta_\e\id_{\gamma_1}(x)+(1-\delta_\e) F_\e(x)$, $x\in\R$, which is rewritten in the form
\begin{eqnarray}\label{eq_Fecirc}
	F_\e^\circ(t)=  (\delta_\e+(1-\delta_\e)c_d)\id_{\gamma_1}(x) +(1-\delta_\e)(1-c_d)F_{a,\e}(x),\quad x\in\R.
\end{eqnarray}
It is seen that $F_\e^\circ\in \DFdisc(L)\oplus\DFabs(S)$. By Fact~\ref{fact_discLabs}, where we set $F_d:=\id_{\gamma_1}$ (i.e. $f_d(t)=e^{it\gamma_1}\ne 0$, $t\in\R$), we have $F_\e^\circ\in\RID$ and, consequently, $F_\e^\circ\in\RIDmixdabs(L,S)$. Let us estimate the variation of  $F-F_\e^\circ$. Observe that
\begin{eqnarray}\label{ineq_FFecirc}
	\|F-F_\e^\circ\|\leqslant \|F-F_\e\|+\|F_\e-F_\e^\circ\|.
\end{eqnarray}
For the first term of the right-hand side we have
\begin{eqnarray}
	\|F-F_\e\|&=&(1-c_d)\|F_a-F_{a,\e}\| \leqslant  \|F_a-F_{a,\e}\|\nonumber\\
	&=&\|F_a-c_\e F_{a,\e}-(1-c_\e) F_{a,\e}\|\nonumber\\
	&\leqslant& \|F_a-c_\e F_{a,\e}\|+(1-c_\e) \|F_{a,\e}\|\nonumber\\
	&=& \|F_a-c_\e F_{a,\e}\|+1-c_\e.\label{ineq_FFe1a}
\end{eqnarray}
Here  $1-c_\e<\e$ and $\|F_a-c_\e F_{a,\e}\|<\e$ as in the proof of Theorem~\ref{th_AbsTV} (see formulas \eqref{ineq_int_p_e} and \eqref{ineq_FceFe}). Therefore $\|F-F_\e\|<2\e$. We next estimate the second term in the right-hand side of \eqref{ineq_FFecirc}:
\begin{eqnarray}\label{ineq_FFecirc1a}
	\|F_\e-F_\e^\circ\|=\|F_\e - \delta_\e\id_{\gamma_1}-(1-\delta_\e) F_\e\|=\delta_\e\| F_\e- \id_{\gamma_1}\|\leqslant 2\delta_\e<2\e.
\end{eqnarray}
Thus we obtain $\|F-F_\e^\circ\|\leqslant 4\e$.

\textbf{1.b)} We now assume that $\cext(F_\e)=\gamma_1$. Since $\Supp(F_\e)=\{\gamma_1\}\cup \Supp(F_{a,\e})$ and $\lext(F_{a,\e})<\rext(F_{a,\e})$, the inequalities $\gamma_1\leqslant\lext(F_{a,\e})$ and $\gamma_1\geqslant\rext(F_{a,\e})$ can not be valid. Then
\begin{eqnarray}\label{ineq_gamma1}
	\lext(F_{a,\e})< \gamma_1<\rext(F_{a,\e})
\end{eqnarray}
and hence $\gamma_1=\cext(F_{a,\e})$. In order to apply Lemma~\ref{lm_key} here, we will slightly change the function $F_{a,\e}$. Namely, we choose a real number $\hat{r}_\e$  such that $\gamma_1<\hat{r}_\e<\rext(F_{a,\e})$ and $\int_{[\hat{r}_\e,r_\e]} p(x)\dd x<\e$. Here we define $\hat{c}_\e:=\int_{[-r_\e,\hat{r}_\e]} p(x)\dd x$ and   $\hat{p}_\e(x):= \tfrac{1}{\hat{c}_\e}\,p(x)\iid_{[-r_\e,\hat{r}_\e]}(x)$, $x\in\R$. Since $\e<\tfrac{1}{2}$, the constant $\hat{c}_\e$ is (strictly) positive:
\begin{eqnarray}\label{ineq_hatce}
	\hat{c}_\e=\int_{[-r_\e,r_\e]} p(x)\dd x-\int_{[\hat{r}_\e,r_\e]} p(x)\dd x > (1-\e)-\e=1-2\e>0.
\end{eqnarray}
So $\hat{p}_\e$ is the density function of a probability law. We denote by 
$\hat{F}_{a,\e}$ and $\hat{f}_{a,\e}$ the distribution function and the characteristic function corresponding to $\hat{p}_\e$, respectively. Observe that  $\Supp(\hat{F}_{a,\e})\subset\Supp(F_{a,\e})\subset\Supp(F_a)\subset S$ and $\lext(\hat{F}_{a,\e})=\lext(F_{a,\e})$, however $\rext(\hat{F}_{a,\e})\leqslant \hat{r}_\e<\rext(F_{a,\e})$. Next, we set 
\begin{eqnarray*}
	\hat{F}_\e(x):=c_d \id_{\gamma_1}(x)+(1-c_d) \hat{F}_{a,\e}(x),\quad x\in\R.
\end{eqnarray*}
Here  $\Supp(\hat{F}_\e)=\{\gamma_1\} \cup \Supp(\hat{F}_{a,\e})$ and, on account of \eqref{ineq_gamma1}, we have
\begin{eqnarray*}
	\cext(\hat{F}_\e)&=&\dfrac{1}{2}\bigl(\lext(F_{a,\e})+\max\bigl\{\rext(\hat{F}_{a,\e}),\gamma_1\bigr\}\bigr)\\
	&<&\dfrac{1}{2}\bigl(\lext(F_{a,\e})+\rext(F_{a,\e})\bigr)=\cext(F_{a,\e})=\gamma_1.
\end{eqnarray*}
Thus $\cext(\hat{F}_\e)\ne\gamma_1$ and we come to the case \textbf{1.a)}. Similarly to that case, we define a new distribution function $\hat{F}_\e^\circ(x):=\delta_\e\id_{\gamma_1}(x)+(1-\delta_\e) \hat{F}_\e(x)$, $x\in\R$, whose characteristic function has no zeroes on the real line for some $\delta_\e\in(0,\e)$. So $\hat{F}_\e^\circ\in\RIDmixdabs(L,S)$ too. The estimates for the variation of $F-\hat{F}_\e^\circ$ are analogous:
\begin{eqnarray*}
	\|F-\hat{F}_\e^\circ\|\leqslant \|F-\hat{F}_\e\|+\|\hat{F}_\e-\hat{F}_\e^\circ\|,
\end{eqnarray*}
where $\|F-\hat{F}_\e\|\leqslant  \|F_a-\hat{c}_\e \hat{F}_{a,\e}\|+1-\hat{c}_\e<4\e$ (see \eqref{ineq_FceFe}, \eqref{ineq_FFe1a} and \eqref{ineq_hatce}), and $\|\hat{F}_\e-\hat{F}_\e^\circ\|=\delta_\e\| \hat{F}_\e- \id_{\gamma_1}\|\leqslant 2\delta_\e<2\e$ (see \eqref{ineq_FFecirc1a}). So we get $\|F-\hat{F}_\e^\circ\|\leqslant 6\e$.

\textbf{2)} We now turn to the case, when  $L$ contains at least two distinct points $a+bk_1$ and $a+bk_2$ with integers $k_1$ and $k_2$ (we set $k_1<k_2$). Since we have already handled  the case of singleton set $\Supp(F_d)$,  we can think that $\Supp(F_d)$ contains these points too without loss of generality. According to Theorem~\ref{th_DiscLTV}, for any $\e\in(0,\tfrac{1}{2})$ we can find   $F_{d,\e}\in\RIDdisc(L)$ such that $\|F_d-F_{d,\e}\|<\e$. We choose such an $F_{d,\e}$ with the bounded set $\Supp(F_{d,\e})$ (see the proof of Theorem~\ref{th_DiscLTV}). Let $f_{d,\e}$ be the characteristic function of $F_{d,\e}$. Due to Fact~\ref{fact_disc}, $\mu_{d,\e}:= \inf_{t\in\R}|f_{d,\e}|>0$. Next, we define $F_{a,\e}$ and $f_{a,\e}$ as above, see \textbf{1)}. Recall that $\Supp(F_{a,\e})$ is bounded and it is contained in $S$.  We also have the estimate $\|F_a-F_{a,\e}\|<2\e$, see \eqref{ineq_FFe1a} in \textbf{1.a)}.  Let us define the distribution function
\begin{eqnarray}\label{def_Fe2}
	F_\e(x):=c_d F_{d,\e}(x)+(1-c_d)F_{a,\e}(x),\quad x\in\R,
\end{eqnarray}
with  the corresponding characteristic function $f_\e$. It is clear that $\Supp(F_\e)$ is bounded. Since $L$ contains at least two points, there is $\gamma_2\in L$ such that $\gamma_2\ne \cext(F_\e)$. By Lemma~\ref{lm_key}, for any $\tau\in(0,1]$  there exists $\delta\in(0,\tau)$ such that
\begin{eqnarray}\label{ineq_fecirc}
	f^\circ_\e(t):=\delta e^{it\gamma_2}+(1-\delta) f_\e(t)\ne 0\quad\text{for any}\quad t\in\R.
\end{eqnarray}
For any fixed $\delta$ the function $f^\circ_\e$ is characteristic for the distribution function $F_\e^\circ(x):=\delta \id_{\gamma_2}(x)+(1-\delta) F_\e(x)$, $x\in\R$. Using \eqref{def_Fe2}, we have
\begin{eqnarray*}
	F_\e^\circ(x)&=&\delta \id_{\gamma_2}(x)+(1-\delta)c_d F_{d,\e}(x)+(1-\delta)(1-c_d)F_{a,\e}(x)\\
	&=&(\delta+c_d-\delta c_d) F_{d,\e}^\circ(x)+(1-\delta-c_d+\delta c_d)F_{a,\e}(x),\quad x\in\R,
\end{eqnarray*} 
where 
\begin{eqnarray*}
	F_{d,\e}^\circ(x):=\dfrac{\delta \id_{\gamma_2}(x)+(1-\delta)c_d F_{d,\e}(x)}{\delta+c_d-\delta c_d},\quad x\in\R.
\end{eqnarray*}
It is easily seen that $F_\e^\circ\in \DFdisc(L)\oplus\DFabs(S)$. If $c_d=0$, i.e. $F_\e^\circ(x)=\delta \id_{\gamma_2}(x)+(1-\delta)F_{a,\e}(x)$, $x\in\R$, then we choose $\tau:=\e$ and so $\delta=\delta_\e<\e$. This case has already been covered  in \textbf{1.a)} (see \eqref{eq_Fecirc} with $c_d=0$). So we have $F_\e^\circ\in\RIDmixdabs(L,S)$ and  $\|F-F_\e^\circ\|\leqslant 4\e$. Next, if $c_d>0$ then we choose $\delta=\delta_\e<\tau:= \min\{\e,(1-\e)c_d \mu_{d,\e}\}$ and consider the characteristic function of $F_{d,\e}^\circ$:
\begin{eqnarray*}
	f_{d,\e}^\circ(t):=\dfrac{\delta_\e e^{it\gamma_2}+(1-\delta_\e)c_d f_{d,\e}(t)}{\delta_\e+c_d-\delta_\e c_d},\quad t\in\R.
\end{eqnarray*}
Observe that
\begin{eqnarray*}
	|f_{d,\e}^\circ(t)|\geqslant \dfrac{(1-\delta_\e)c_d |f_{d,\e}(t)|-\delta_\e}{\delta_\e+c_d-\delta_\e c_d}\geqslant\dfrac{(1-\e)c_d \mu_{d,\e}-\delta_\e}{\delta_\e+c_d-\delta_\e c_d}\geqslant\dfrac{\tau-\delta_\e}{\delta_\e+c_d-\delta_\e c_d}>0.
\end{eqnarray*}
From this and \eqref{ineq_fecirc} we conclude that $F_\e^\circ\in\RID$ due to Fact~\ref{fact_discLabs}. Hence $F_\e^\circ\in\RIDmixdabs(L,S)$. Next, we have the following estimates for the variations:
\begin{eqnarray*}
	\|F-F_\e^\circ\|\leqslant \|F-F_\e\|+\|F_\e-F_\e^\circ\|,
\end{eqnarray*}
where
\begin{eqnarray*}
	\|F-F_\e\|&=&  \|c_d (F_d-F_{d,\e})+(1-c_d)(F_a-F_{a,\e})\|\\
	&\leqslant& c_d\|F_d-F_{d,\e}\| +(1-c_d) \|F_a-F_{a,\e}\|\\
	&<& c_d \e+(1-c_d) 2\e< 2\e,
\end{eqnarray*}
and $\|F_\e-F_\e^\circ\|=\delta_\e\| F_\e- \id_{\gamma_2}\|\leqslant 2\delta_\e<2\e$ (see \eqref{ineq_FFecirc1a}). Therefore $\|F-F_\e^\circ\|<4\e$.

Thus for any $F\in \DFdisc(L)\oplus\DFabs(S)$ and for any case, \textbf{1)} or \textbf{2)},  we can find $F_\e^\circ\in\RIDmixdabs(L,S)$, which is arbitrarily close to $F$ in variation.  So we come the required denseness of $\RIDmixdabs(L,S)$ in $\DFdisc(L)\oplus\DFabs(S)$.\quad $\Box$\\

Thus, according to Theorems~\ref{th_AbsTV}--\ref{th_Mixdabs}, discrete lattice laws,  absolutely continuous laws and their mixtures (say the distributions, which are the most important for practice) can be  approximated arbitararily well by rational-infinitely divisible distributions (of the same types, correspondingly)  with respect to very strong distance in total variation, i.e. a distribution with one of these types is almost indistinguishable with some rational-infinitely law of this type.

\section{The ``negative'' results}

Here we follow the notation from the previous sections. We focus on the class $\DFdiscNL$ of  discrete non-lattice distribution functions.

\begin{Theorem}\label{th_DiscreteTV}
	The class $\RIDdisc$ is not dense in variation in the class $\DFdiscNL$.
\end{Theorem}
The following obvious corollary is opposite to Corollary~\ref{co_th_DiscLTV} for the lattice case.
\begin{Corollary}
	The class $\RIDdiscNL$ is not dense in variation in the class $\DFdiscNL$.
\end{Corollary}

\noindent\textbf{Proof of Theorem~\ref{th_DiscreteTV}.}
Let $\alpha$ be a positive irrational number and
\begin{eqnarray*}
	F(x):=\dfrac{1}{3}\,\id_1(x)+\dfrac{1}{3}\,\id_{\alpha}(x)+\dfrac{1}{3}\,\id_{1+\alpha}(x),\quad x\in\R.
\end{eqnarray*}
It is easy to see that $F\in\DFdiscNL$. We assume that there is a sequence $(F_n)_{n\in\N}$ from $\RIDdisc$ such that 
\begin{eqnarray}\label{conv_FnF}
	\|F_n-F\|\to 0\quad\text{as}\quad n\to\infty.
\end{eqnarray}
We represent every $F_n$ in the following form:
\begin{eqnarray*}
	F_n(x)=p_{n,1}\id_{1}(x)+p_{n,\alpha}\id_\alpha(x)+p_{n,1+\alpha} \id_{1+\alpha}(x) +\sum_{k=1}^\infty p_{n,u_k} \id_{u_k}(x),\quad x\in\R,\quad n\in\N,
\end{eqnarray*}
where for every $n\in\N$ we set $p_{n,u_k}\geqslant 0$, $k\in\N$, and $p_{n,1}+p_{n,\alpha} +p_{n,1+\alpha}+\sum_{k=1}^\infty p_{n,u_k}=1$. Here, without loss of generality, $\U:=\{u_1, u_2, \ldots, u_k, \ldots\}$ is an infinite countable set of distinct real numbers such that $\{1, \alpha, 1+\alpha\}\cap\U=\varnothing$ and the set $\{1, \alpha, 1+ \alpha\}\cup \U$ contains all the sets $\Supp(F_n)$, $n\in\N$. So we have
\begin{eqnarray}\label{eq_FnF}
	\|F_n-F\|=\bigl|p_{n,1}-1/3\bigr|+\bigl|\,p_{n,\alpha} -1/3\bigr|+\bigl|p_{n,1+\alpha}-1/3\bigr| +\sum_{k=1}^\infty p_{n,u_k},\quad n\in\N.
\end{eqnarray}

Let $f_n$ denote the characteristic function of $F_n$ for every $n\in\N$. Then 
\begin{eqnarray*}
	f_n(t)=p_{n,1}e^{it}+p_{n,\alpha}e^{it\alpha}+p_{n,1+\alpha} e^{it(1+\alpha)} +\sum_{k=1}^\infty p_{n,u_k} e^{itu_k},\quad t\in\R,\quad n\in\N.
\end{eqnarray*}
Since $F_n\in\RID_d$, we know that $\mu_n:=\inf_{t\in\R}|f_n(t)|>0$, $n\in\N$, by Fact~\ref{fact_disc}. Next, for every $n\in\N$ we choose $K_n\in\N$ such that $\sum_{k=K_n+1}^\infty p_{n,u_k}<\mu_n/2$, and we set
\begin{eqnarray*}
	\tilde{f}_n(t):=p_{n,1}e^{it}+p_{n,\alpha}e^{it\alpha}+p_{n,1+\alpha} e^{it(1+\alpha)} +\sum_{k=1}^{K_n} p_{n,u_k} e^{itu_k},\quad t\in\R,\quad n\in\N.
\end{eqnarray*}  
Observe that
\begin{eqnarray*}
	\inf_{t\in\R}|\tilde{f}_n(t)|\geqslant \inf_{t\in\R}|f_n(t)|-\sup_{t\in\R}|f_n(t)-\tilde{f}_n(t)|
	=\mu_n-\sup_{t\in\R}\biggl|\sum_{k=K_n+1}^\infty p_{n,u_k}e^{itu_k}\biggr|
\end{eqnarray*}
and hence we have
\begin{eqnarray}\label{ineq_inftildefn}
	\inf_{t\in\R}|\tilde{f}_n(t)|\geqslant \mu_n-\sum_{k=K_n+1}^\infty p_{n,u_k}>\mu_n-\dfrac{\mu_n}{2}=\dfrac{\mu_n}{2}>0.
\end{eqnarray}

We next select a basis of the system $\U_0:=\{1, \alpha, 1+\alpha\}\cup\,\U$ over $\Q$ by the following procedure. We set $\beta_1:=1$ and $\beta_2:=\alpha$. We next sequentially  cross out all numbers $y$ of $\U_0$ satisfying the relation $r_1 \beta_1+r_2\beta_2+r_3 y=0$ with any $r_1$, $r_2$ and $r_3\ne 0$ from $\Q$.   Of course, $1$, $\alpha$ and $1+\alpha$ will be cross out. Let $\beta_3$ be the first term of the sequence $(u_k)_{k\in\N}$ ($\U$ with order), which was not crossed out.  Next, we cross out from the set of the numbers, which have not been crossed out yet, all numbers $y$ satisfying the relation $r_1 \beta_1+r_2\beta_2+r_3 \beta_3 +r_4 y=0$ with any $r_1$, $r_2$, $r_3$ and $r_4\ne 0$ from $\Q$. We denote by $\beta_4$ the first number in $(u_k)_{k\in\N}$, which was not crossed out. By continuing this process analogously, we get the collection $\beta_1$, $\beta_2$, \ldots,  $\beta_k$,\ldots, which may be finite (if all elements of $\U_0$ will be crossed out at some step of the process) or infinite. By construction,  it is seen that $\beta_1$, $\beta_2$, \ldots are non-zero and they linearly independent over $\Q$, i.e. if $r_1\beta_1+r_2\beta_2+\ldots+r_m\beta_m=0$ with $m\in\N$ and $r_1, r_2, \ldots, r_m\in \Q$, then $r_1=r_2=\ldots=r_m=0$. In addition, every $u_k\in\U_0$ is represented by the finite linear combination
\begin{eqnarray*}
	u_k= a_{k,1}\beta_1+a_{k,2}\beta_2+\ldots+a_{k,m_k}\beta_{m_k},\quad k\in\N,
\end{eqnarray*}
with some $m_k\in\N$ and $a_{k,1}, a_{k,2}, \ldots, a_{k,m_k}\in \Q$. Of course, this representation is unique (if the order of summation and the summands with zero coefficients $a_{k,j}$ are disregarded). Since zero values are admissible for $a_{k,j}$, we will assume that $m_k\geqslant 2$ for convenience.

We now return to the functions $\tilde{f}_n$, $n\in\N$. Let us fix $n\in\N$ and find the minimal positive integer $\varkappa_{n}$ such that $A_{k,j}^n:=\varkappa_{n} a_{k,j}\in \Z$ for any $k=1,\ldots, K_n$ and $j=1,\ldots, m_k$. We set $\bar\beta_j^n:=\beta_j/\varkappa_{n}$ for every admissible $j$. Therefore we write
\begin{eqnarray*}
	u_k= A_{k,1}^n\bar\beta_1^n+A_{k,2}^n\bar\beta_2^n+\ldots+A_{k,m_k}^n\bar\beta_{m_k}^n,\quad k=1,\ldots, K_n,
\end{eqnarray*}
Thus we can represent $\tilde{f}_n$ in the following form
\begin{eqnarray*}
	\tilde{f}_n(t)&=&p_{n,1}\exp\bigl\{it\varkappa_{n}\bar\beta_1^n\bigr\}+p_{n,\alpha}\exp\bigl\{it\varkappa_{n}\bar\beta_2^n\bigr\}+p_{n,1+\alpha} \exp\bigl\{it\bigl(\varkappa_{n}\bar\beta_1^n+\varkappa_{n}\bar\beta_2^n\bigr)\bigr\}\\
	&&{} +\sum_{k=1}^{K_n} p_{n,u_k} \exp\biggl\{it\sum_{j=1}^{m_k} A_{k,j}^n\bar\beta_j^n\biggr\},\quad t\in\R.
\end{eqnarray*}
Let $M_n:=\max\{m_k: k=1,\ldots, K_n\}$. We introduce the function
\begin{eqnarray*}
	\tilde{\varphi}_n(t_1,\ldots,t_{M_n})&=&p_{n,1}\exp\bigl\{it_1\varkappa_{n}\bar\beta_1^n\bigr\}+p_{n,\alpha}\exp\bigl\{it_2\varkappa_{n}\bar\beta_2^n\bigr\}+p_{n,1+\alpha} \exp\bigl\{i\bigl(t_1\varkappa_{n}\bar\beta_1^n+t_2\varkappa_{n}\bar\beta_2^n\bigr)\bigr\}\\
	&&{}+\sum_{k=1}^{K_n} p_{n,u_k} \exp\biggl\{i\sum_{j=1}^{m_k} t_jA_{k,j}^n\bar\beta_j^n\biggr\},\quad t_1,\ldots,t_{M_n}\in\R,
\end{eqnarray*}
for which the function $\tilde{f}_n$ is diagonal, i.e. $\tilde{\varphi}_n(t,\ldots,t)=\tilde{f}_n(t)$, $t\in\R$. Note that $\tilde{\varphi}_n$ is a continuous periodic function with periods $2\pi/\bar\beta_1^n$, . . . ,$2\pi/\bar\beta_{M_n}^n$ over variables
$t_1,\ldots,t_{M_n}$, respectively (recall that all $A_{k,j}^n$ are integers). Since the reciprocals of the periods are linear independent over $\Q$,  the set of values of  $\tilde{f}_n$ is dense in the set of values of $\tilde{\varphi}_n$ by the well known fact from the general theory of almost periodic functions (see \cite{Levitan}, Theorem~2.4.1., p.~116). Hence $|\tilde{\varphi}_n(t_1,\ldots,t_{M_n})|\geqslant \mu_n/2$ for any $t_1,\ldots,t_{M_n}\in\R$ due to \eqref{ineq_inftildefn}. Next, we consider the function
\begin{eqnarray*}
	\tilde{\varphi}_{2,n}(t_1,t_2)&:=&\tilde{\varphi}_n(t_1, t_2/\alpha,0,\ldots,0)\\
	&=&p_{n,1}\exp\bigl\{it_1\varkappa_{n}\bar\beta_1^n\bigr\}+p_{n,\alpha}\exp\bigl\{it_2\varkappa_{n}\bar\beta_2^n/\alpha\bigr\}+p_{n,1+\alpha} \exp\Bigl\{i\bigl(t_1\varkappa_{n}\bar\beta_1^n+t_2\varkappa_{n}\bar\beta_2^n/\alpha\bigr)\Bigr\}\\
	&&{}\quad +\sum_{k=1}^{K_n} p_{n,u_k} \exp\Bigl\{i\bigl(t_1A_{k,1}^n\bar\beta_1^n+t_2A_{k,2}^n\bar\beta_2^n/\alpha \bigr) \Bigr\}\\
	&=& p_{n,1}e^{it_1}+p_{n,\alpha} e^{it_2}+p_{n,1+\alpha}e^{i(t_1+t_2)}+\sum_{k=1}^{K_n} p_{n,u_k} \exp\Bigl\{i\bigl(t_1 a_{k,1}+t_2a_{k,2}\bigr)\Bigr\},\quad t_1,t_2\in\R.
\end{eqnarray*}
By the above, in particular, we have
\begin{eqnarray}\label{ineq_phi2nmun}
	|\tilde{\varphi}_{2,n}(t_1,t_2)|\geqslant \mu_n/2>0\quad\text{for any} \quad t_1,t_2\in\R.
\end{eqnarray}

We now vary index $n$. Let us introduce the function 
\begin{eqnarray*}
	\varphi_{2}(t_1,t_2):=\dfrac{1}{3}\,e^{it_1}+\dfrac{1}{3}\,e^{it_2}+\dfrac{1}{3}\,e^{i(t_1+t_2)},\quad t_1,t_2\in\R,
\end{eqnarray*}
and observe that for any $n\in\N$
\begin{eqnarray*}
	\Delta_n&:=&\sup_{t_1,t_2\in\R}\bigl|\tilde{\varphi}_{2,n}(t_1,t_2) -\varphi_{2}(t_1,t_2) \bigr|\\
	&=&\sup_{t_1,t_2\in\R} \biggl| \bigl(p_{n,1}-1/3\bigr)e^{it_1} +\bigl(p_{n,\alpha}-1/3\bigr)e^{it_2}+\bigl(p_{n,1+\alpha}-1/3\bigr)e^{i(t_1+t_2)}\\
	&&{}+ \sum_{k=1}^{K_n} p_{n,u_k} \exp\Bigl\{i\bigl(t_1 a_{k,1}+t_2a_{k,2}\bigr) \Bigr\}\biggr|\\
	&\leqslant& \bigl|p_{n,1}-1/3\bigr|+\bigl|p_{n,\alpha} -1/3\bigr|+\bigl|p_{n,1+\alpha}-1/3\bigr| +\sum_{k=1}^{K_n} p_{n,u_k}.
\end{eqnarray*}
According to \eqref{conv_FnF} and \eqref{eq_FnF}, we have $\Delta_n\leqslant \|F_n-F\|$ and hence $\Delta_n\to 0$ as $n\to\infty$. Due to Lemma~\ref{lm_Kutlu}, we conclude that $\tilde{\varphi}_{2,n}$ have zeroes in $\R^2$ for all large enough $n\in\N$. However, it obviously contradicts to \eqref{ineq_phi2nmun}.

Thus it is impossible to select a sequence $(F_n)_{n\in\N}$ from $\RID_d$, which converges to $F$ in variation.\quad $\Box$\\

The next theorem is rather non-obvious corollary of the previous one. It is obtained due to the fact from Lemma~\ref{lm_fdsep0} (see the recent results from \cite{KhartovNecess}). 

\begin{Theorem}\label{th_QFdNL}
	The class $\RID$ is not dense in variation in the class $\DFdiscNL$.
\end{Theorem}
\noindent\textbf{Proof of Theorem~\ref{th_QFdNL}.}
Suppose, contrary to our claim, that $\RID$ is dense in variation in $\DFdiscNL$. Then for any $F\in \DFdiscNL$ there exists a sequence $(F_n)_{n\in\N}$ from $\RID$ such that 
\begin{eqnarray}\label{conv_FnF2}
	\|F_n-F\|\to 0\quad\text{as}\quad n\to\infty.
\end{eqnarray}
Let us represent each $F_n$ as the sum of its discrete and continuous parts:
\begin{eqnarray*}
	F_n(x)=c_{d,n} F_{d,n}(x)+(1-c_{d,n}) F_{c,n}(x),\quad x\in\R,\quad n\in\N,
\end{eqnarray*}
where $F_{d,n}$ is a discrete distribution function, $F_{c,n}$ is a continuous one, and $c_{d,n}\in[0,1]$. Then
\begin{eqnarray*}
	\|F_n-F\|=\|c_{d,n} F_{d,n}-F\|+\|(1-c_{d,n}) F_{c,n} \|=\|c_{d,n} F_{d,n}-F\|+1-c_{d,n},\quad n\in\N.
\end{eqnarray*}
From \eqref{conv_FnF2} we have $\|c_{d,n} F_{d,n}-F\|\to 0$ and $c_{d,n}\to 1$ as $n\to\infty$. It implies that $\|F_{d,n}-F\|\to 0$ as $n\to\infty$. Indeed, it follows from the estimates:
\begin{eqnarray*}
	\|F_{d,n}-F\|\leqslant \|(1-c_{d,n})F_{d,n}\|+\|c_{d,n}F_{d,n}-F\|=1-c_{d,n}+\|c_{d,n}F_{d,n}-F\|,\quad n\in\N. 
\end{eqnarray*}
Next, we know that $F_n\in\RID$ and $c_{d,n}>0$ for all large enough $n\in\N$. Hence, according to Lemma~\ref{lm_fdsep0}, the characteristic function  of each $F_{d,n}$, say $f_{d,n}$, is separated from zero, i.e. $\inf_{t\in\R}|f_{d,n}(t)|>0$. By Fact~\ref{fact_disc}, it means that $F_{d,n}\in\RID_d$  for every $n\in\N$. 

Thus for arbitrary $F\in \DFdiscNL$ we can always find the sequence $(F_{d,n})_{n\in\N}$ from $\RID_d$, which converges in variation to $F$. But it is not true by Theorem~\ref{th_DiscreteTV}. Hence the initial assumption is false, i.e. $\RID$ is not dense in variation in $\DFdiscNL$.\quad $\Box$\\

Theorem~\ref{th_QFdNL} immediately implies the following important assertion, which we formulate as theorem.

\begin{Theorem}\label{th_QF}
	The class $\RID$ is not dense in variation in the class $\DF$.
\end{Theorem}

Due to the results by Kutlu (see \cite{Kutlu}), this theorem is also valid for the multivariate analog of $\RID$, because, for instance, the multivariate version of Corollary~\ref{co_th_DiscLTV} is false. However, does Theorem~\ref{th_AbsTV} (or Corollary~\ref{co_th_AbsTV}) hold in the multivariate case? This is an interesting question.

\end{document}